\documentclass[a4paper]{article}
\usepackage[utf8]{inputenc}
\usepackage{amssymb,amsbsy,amsmath,amsfonts,amscd}
\usepackage{euscript}
\usepackage{a4wide}
\usepackage{graphicx}
\usepackage[usenames,dvipsnames]{color}
\usepackage{comment}

\title{Some remarks on simplified double porosity model of immiscible incompressible two-phase flow}
\author{M. Jurak$^1$, L. Pankratov$^2$, A. Vrba\v ski$^3$}
\date{\today}

\begin{document}
	
	\maketitle
	
	
	%

	\def \ve{\varepsilon} 
	\def \eps{\epsilon}
	\def \del{\delta}
	\def \ved{{\varepsilon,\del}}%
	\def \KSD{\K^{\star, \delta}}
	\def \lm{\lambda}
	\def \Lm{\Lambda}
	\def \limsup{\mathop{\overline {\rm lim}}}
	\def \liminf{\mathop{\underline{\rm lim}}}
	\def \la{\left\langle\rule{0pt}{3em}}
	\def \ra{\right\rangle}
	\newcommand{\Ab}[1]{(A.#1)}
	\newcommand{\Hb}[1]{(\textsf{A.#1})}
	\def \gr{\nabla}
	\def \pt{\partial}
	\def \ptt{\partial_t}
	\def \dv{\mbox{$\mbox{div}$}}
	\def \div{{\rm div}\,}
	\def \to{\rightarrow}
	\def \tow{\rightharpoonup}
	\def \eqdef{\stackrel {\rm def} {=}}
	\def \be{\begin{equation}}
		\def \ee{\end{equation}}
	\def \ds{\displaystyle}
	\def \meas{{\rm meas}\,}
	\def \R{{\mathbb R}} \def \C{{\Bbb C}} \def \Z{{\Bbb Z}}
	\def \I{{\Bbb I}} \def \N{{\Bbb N}} \def \K{{\mathbb K}}
	\def \bs{\boldsymbol}
	\def \ww{\mathfrak{w}}
	\def \rr{\mathfrak{r}}
	\def \ss{\mathfrak{S}}
	\def \oz{\overline{z}}
	\def \ox{\overline{x}}
	\def \mx{\mathsf m}
	\def \fr{\mathsf f}
	\def \mm{\mathfrak{M}}
	\def \ff{\mathfrak{F}}
	\def \ge{\mathfrak{g}}
	\def \aa{\mathfrak{A}}
	\def \e1{\mathfrak{e}}
	
	\newtheorem{definition}{Definition}
	\newtheorem{lemma}{Lemma}
	
	\newcommand{\red}[1]{\textcolor{red}{#1}}
	\newcommand{\blue}[1]{\textcolor{blue}{#1}}
	\def\green{\textcolor{green}}
	\def\magenta{\textcolor{magenta}}
	\def\cyan{\textcolor{cyan}}

	\title{}

	\date{}
	\maketitle
	
	\begin{small}
	
	\noindent $^1$ Faculty of Science, University of Zagreb, Bijeni\v cka 30, 10000 Zagreb, Croatia.
	E-mail: {\tt jurak@math.hr}
	
	\noindent $^2$ Laboratory of Fluid Dynamics and Seismic ({\sl RAEP 5top100}), Moscow Institute
	of Physics and Technology, 9 Institutskiy per., Dolgoprudny, Moscow Region, 141700, Russian Federation and 
	 Laboratoire de Math\'ematiques et de leurs Applications, CNRS-UMR 5142 Universit\'e de Pau,
	Av. de l'Universit\'e, 64000 Pau, France. 
	E-mail: {\tt leonid.pankratov@univ-pau.fr}
	
	\noindent $^3$  Faculty of Mining, Geology and Petroleum Engineering, University of Zagreb,
	Pierottijeva 6, 10000 Zagreb, Croatia.
	E-mail: {\tt anja.vrbaski@rgn.unizg.hr}
	\end{small}

	\begin{abstract}
		The paper is devoted to the derivation, by linearization, of simplified
		(fully homogenized) homogenized models
		of an immiscible incompressible two-phase flow in double porosity media in the case of
		thin fissures. In a simplified dual porosity model derived previously by the authors the matrix-fracture source term
		is approximated by a convolution type source term. This approach enables to exclude the
		cell problem, in form of the imbibition equation, from the global double porosity model.
		In this paper we propose a new linear version of the imbibition equation which leads
		to a new simplified dual porosity model. 
		We also present numerical simulations which show that the matrix-fracture exchange term based on this new linearization procedure gives a better
		approximation of the exact one than the corresponding exchange term obtained
		earlier by the authors.
		
	\end{abstract}
	
	
	
	
	\section{Introduction}
	\label{sec-intro}
	
	Naturally fractured reservoirs are characterized by a system of fractures existing within a
	background rock matrix. The fracture system has a low storage capacity and a high conductivity,
	while the matrix block system has a conductivity that is low in comparison with that in the fractures.
	The majority of fluid transport will occur along flow paths through the fissure system, and
	the relative volume and storage capacity of the porous matrix is much larger than that of the fissure
	system (type 2 reservoirs in \cite{Nelson-2001}).
	Multiphase flow in subsurface fractured media offers a particular challenge to numerical
	modeling, since both the flow through the fracture network
	and transfer to a relatively stagnant matrix need to be
	modeled.
	The fractures have a very strong influence on flow and transport but the
	discrepancy between the width of the fractures and other dimensions involved makes
	the inclusion of fractures in a numerical model  difficult
	and costly.
	Possible applications of such systems are in improved oil recovery in
	hydrocarbon reservoirs, non-aqueous phase contaminant transport, nuclear waste
	containment etc. (see, {e.g.,} \cite{Bear-1993}).
	
	Mathematical models of multiphase flow in {fractured-porous} media can be categorized
	into two classes: the dual-continuum models and the discrete-fracture network models.
	Discrete fracture modeling represents each fracture and the matrix as a geometrically well-defined
	entities which are represented explicitly in the computational domain.
	This approach leads to the  most accurate and physically realistic model at the expense of highly
	refined or hybrid grid \cite{Firoo-2005, Bastian-06}.
	
	In dual-continuum model one do not consider specific known fractures
	that might be included individually in the model but a network of small interconnected
	fractures with certain degree of regularity.
	There are two large differences in scale in the fractures and the blocks:
	fracture width is very small compared to the scale of the domain and fracture permeability
	is much larger than the permeability of surrounding material.
	In dual-continuum approach a sort of averaging process is applied in order to obtain
	simplified description of matrix-fracture system and its interactions. In the case of stagnant
	flow in the matrix averaging leads to a dual porosity model that was first {obtained} in
	\cite{BZK-1960, WarrenRoot}, where the model was described as a
	phenomenological model deduced experimentally. In the framework of the model presented
	in these papers the fissures are assumed to have a negligible volume with respect
	to the volume of the whole reservoir (see \cite{BZK-1960} {or \cite{RAG} Chapter 5}).
	
	From the mathematical point of view, the usual double-porosity model
	(or $\ve^2$-model) assumes that the width of the fracture containing highly
	permeable porous media is of the same order as the block size.
	The related homogenization problem was first studied in \cite{Arbogas-90}, and was then
	revisited in the mathematical literature by many other authors (see, e.g.
	\cite{bmp,bgpp,brad,cat,Hornung-97,kh,panf,san,ak-lp-AA,ainouz,BA-LP-M2AS,AMP-B},
	and references therein).

	Homogenization of linear and nonlinear parabolic equations, in particular,
	parabolic equations with high contrast coefficients (see the references above) is
	a long-standing question. In this paper we deal with periodic homogenization of
	double-porosity-type problems in a special case of asymptotically small volume
	fraction of highly conductive part of the medium. In other words, when the
	permeability coefficient of the matrix part vanishes everywhere in the corresponding
	domain except a set of asymptotically small measure. The first results on the subject
	were obtained independently and by different approaches in  \cite{bcp} and
	\cite{pr}. In \cite{pr} the highly conductive part of the medium (the fissures system)
	is modelized by the only one small parameter $\ve$ which characterizes the scale of the
	microstructure. The main feature of the homogenized model in this case is that it does not
	involve neither a cell problem for the calculation of the additional source term nor
	a cell problem for the calculation of the global permeability tensor.
	The further results on the subject are connected with the
	application of the method of two small parameters $\ve,\delta$ proposed in \cite{bp,cior}.
	In this case the parameter $\ve$ describes the periodicity of the fractured-porous medium 
	and $\delta$ is the relative thickness or opening of the fracture. The homogenization
	process then splits in two steps. First, one shows that the corresponding problem admits
	homogenization as $\ve \to 0$ and then in the second step it is necessary to pass to the
	limit as $\delta \to 0$. Thus the effective system does not depend on $\ve,\delta$.
	We refer here to \cite{rthin-2,th-layer-M3AS} and the references therein.
	The first result on the homogenization of the two-phase flow in double porosity media
	was obtained in \cite{JPV} by the method of two small parameters. 
	In the global double porosity $\delta$-problem obtained in \cite{JPV} after passing to the limit as $\ve\to0$, additional matrix-fracture source terms are present that are given implicitly via solutions of a nonlinear local boundary value problem known as the \textit{imbibition equation}. The nonlinearity of the imbibition equation causes difficulties in numerical simulations of the model since no analytic expression of the matrix-fracture source term is available. In order to overcome this issue, one can linearize the imbibition equation and then express the matrix-fracture source term explicitly from the linearized equation. In \cite{JPV} the imbibition equation is linearized by using an appropriate constant, as suggested by \cite{Arb-simpl}. In this paper we present a new, variable and more general linearization of the imbibition equation which gives a new simplified dual porosity model. We also display numerical simulations comparing the matrix-fracture exchange
	term calculated by solving nonlinear imbibition equation to the matrix-fracture exchange terms given by two different linearization procedures.

	The rest of the paper is organized as follows. In Section \ref{dual-po-model} we present
	a dual porosity model of the two-phase flow in the case of thin fissures, where the opening
	of the fissure is described by a small parameter $\delta$. Section \ref{linear-sect} is
	devoted to the linearization of the imbibition equation which is involved in the global
	dual porosity model as a local problem in the matrix block. In Section \ref{sec:decoupling}
	we present the simplified or fully homogenized dual porosity models.
	In Section \ref{sec:num_comp} we present numerical simulations comparing the matrix-fracture exchange
	term calculated by solving nonlinear imbibition equation to the matrix-fracture exchange terms given by different linearization procedures. The discretization of the effective system including the matrix-fracture source terms obtained by constant linearization is proposed in Section \ref{sec:discrete1}, while in Section \ref{sec:discrete2} we present the discretization of the effective system including the matrix-fracture source terms obtained by variable linearization.

	\section{Dual porosity model}
	\label{dual-po-model}

	In this section we present a dual porosity model of incompressible two-phase  flow
	derived rigorously by the homogenization theory in \cite{ADH91,BLM,Yeh06}.
	Namely, we consider the reservoir $\Omega\subset \mathbb{R}^d$ of the characteristic length $L$ 
	composed of the matrix blocks with the characteristic length $l$  and highly permeable network of fractures.
	The block size iz small compared to the size of the flow domain, i.e., $\varepsilon = l/L \ll 1$ is a 
	small parameter which goes to zero. 
	The thickness of the fractures is supposed to be of order $l\delta $, where $0 < \ve \ll \delta < 1$  
	is a second small parameter. The porosities of the blocks and the fractures are supposed to be constant 
	and are denoted by $\Phi_m$ and $\Phi_f$ respectively. The permeabilities of the blocks and the fractures are 
	highly contrasted. In the derivation of the dual porosity model it is assumed that if the fracture porosity is
	$k_f$, then the matrix porosity is $(l\delta)^2 k_m$, where $k_f$ and $k_m$ are of the same order.

	If we neglect gravitational segregation the dual porosity model obtained by homogenization as $\ve\to0$ (see \cite{ADH91,BLM,Yeh06}) can be written in the form (see \cite{JPV}):
	\begin{equation}
		\label{H-1}
		\left\{
		\begin{array}[c]{ll}
			\ds
			\Phi_{*}^\delta \frac{\pt S_f^\delta}{\pt t} -
			{\rm div}\, \bigg( \mathbb{K}_*^\delta \lm_{w,f} (S_f^\delta) \gr P_{w,f}^\delta \bigg) =
			{\EuScript Q}_w^\delta  \\
			\ds
			- \Phi_{*}^\delta \frac{\pt S_f^\delta}{\pt t} -
			{\rm div}\, \bigg(\mathbb{K}_*^\delta \lm_{n,f} (S_f^\delta) \gr P_{n,f}^\delta \bigg) =
			{\EuScript Q}_n^\delta, \\
			P_{c,f}(S_f^\delta) = P_{n,f}^\delta - P_{w,f}^\delta,
		\end{array}
		\right.
	\end{equation}
	where $S_f^\delta$, $P_{w,f}^\delta$ and $P_{n,f}^\delta$ are  wetting phase saturation and pressure and non wetting 
	phase pressure in the fractures, respectively; $P_{c,f}$, $\lm_{w,f}$ and $\lm_{n,f}$ are the capillary pressure function 
	and the phase mobility functions in the  fractures; $\Phi_{*}^\delta$ and 
	$\mathbb{K}_*^\delta$ are the effective porosity and permeability of the matrix-fracture system. 
	The terms ${\EuScript Q}_w^\delta$ and ${\EuScript Q}_n^\delta$ are the wetting phase and the non wetting phase 
	source terms modeling the phase mass transfer  from the matrix to the fracture system governed by the
	capillary imbibition. 
	
	In this paper as in \cite{JPV}, we adopt the Warren-Root idealization of the
	fractured media, i.e., each matrix block is similar to a rectangular one and all
	the blocks are {surrounded} by the fractures.
	{Then} we introduce the reference cell $Y = (0, 1)^d$ which is decomposed in
	matrix and fracture {parts, where the matrix part, $Y_m^{\del}$,} is an open
	cube with edge length $1 - \delta $ {($0 < \delta \ll 1$)},
	and  $Y_f^{\del} = Y\setminus Y_m^{\del}$ represents the fracture part.
	The flow domain $\Omega$ is assumed to be covered by a pavement of cells $l Y$.

	The effective porosity $\Phi_{*}^\delta$ can be expressed as 
	$$
	\Phi_{*}^\delta= \Phi_f\, \frac{\mbox{vol}(Y_f^{\del})}{\mbox{vol}(Y_m^{\del})}.
	$$
	Moreover, $\Phi_{*}^\delta$ is of order $\delta$. The effective permeability tensor can be
	expressed using the solutions of certain {\sl cell problems} (see \cite{JPV} for more details)
	and it is again of order $\delta$.

	The matrix-fracture source terms ${\EuScript Q}_w^{\del}$,  ${\EuScript Q}_n^{\del}$ are given by:
	\begin{equation}
		\label{H-4}
		{\EuScript Q}_w^{\del}(x,t) \eqdef - \frac{\Phi_m}{\mbox{vol}(Y_m^{\del})} \int\limits_{Y_m^{\del}}
		\frac{\pt S_m^{\del}}{\pt t}(x, y, t)\, dy
		= - {\EuScript Q}_n^{\del}(x,t),
	\end{equation}
	where the function $S_m^{\del}(x, y, t)$ is the matrix block saturation defined for each point $x\in\Omega$ as a solution of  the problem
	\begin{equation}
		\label{imb-eqn}
		\left\{
		\begin{array}[c]{ll}
			\ds
			\Phi_m \frac{\pt S_m^{\del}}{\pt t}\, - \del^2\, k_m\, \Delta_y \beta_m(S_m^{\del}) = 0
			\quad {\rm in}\,\,  Y_m^{\del}, \\
			S_m^{\del}(x, y, t) = \mathcal{P}(S_f^{\del}(x,t)) \quad {\rm on}\,\, \Gamma^{\del}, \\
			S_m^{\del}(x, y, 0) = S_m^0(x) \quad {\rm in}\,\, Y_m^{\del},
		\end{array}
		\right.
	\end{equation}
	where $\Gamma^{\del}$ stands for the interface between the matrix and fracture parts of the
	cell $Y$; $ \mathcal{P}(S) \eqdef P_{c,m}^{-1}(P_{c,f}(S))$ and 
	\begin{equation}
		\label{upsi-1}
		\beta_m(s) \eqdef\int\limits_0^s \alpha_m(\xi)\, d\xi,
		\quad \textrm{ where } \, \,
		\alpha_m(s) \eqdef \frac{\lm_{w,m}(s)\, \lm_{n,m}(s)} {\lm_m(s)} | P^\prime_{c,m}(s) |.
	\end{equation}
	Here $\lm_m(s) = \lm_{w,m}(s) + \lm_{n,m}(s)$. Equation $\eqref{imb-eqn}_1$ is known as the {\em imbibition equation}.
	
	We note that the matrix-fracture source terms ${\EuScript Q}_w^{\del}, {\EuScript Q}_n^{\del}$
	given by (\ref{H-4}) can also be calculated as: 
	\begin{equation*}
		{\EuScript Q}_w^{\del}(x,t) 
		=-\frac{\del^2\, k_m}{\mbox{vol}(Y_m^{\del})} \alpha_m(\mathcal{P}(S_f^{\del}(x,t)) )
		\int\limits_{\partial Y_m^{\del}} \nabla S_m^{\del}\cdot{\bf n}\, dy= - {\EuScript Q}_n^{\del}(x,t).
	\end{equation*}

	\section{Linearization of the imbibition equation } 
\label{linear-sect}

The main difficulty in application of the double porosity model for the two-phase flow is
the fact that the imbibition equation is nonlinear. The nonlinearity does not allow
an analytic solution of the equation and an analytic expression of
{the} corresponding source term. This implies that in the numerical simulation
of the double porosity model  we have to solve the imbibition equation many
times.
{This} problem is especially difficult in the case of
thin fractures, where the solution of the imbibition equation is of {a} boundary
layer type and {a} very refined mesh is needed to resolve it.
{In order to overcome this difficulty, it is possible to linearize} the imbibition equation
and {then use the linearized} equation to express the matrix-fracture transfer term
by an analytic expression.

In this section we explore different ways of linearization of the imbibition equation. The goal is
to derive a linear equation with constant coefficients that can be solved analytically. 

The simplest form of linearization consists in replacing of
a nonlinear function  $\alpha_m(S_m^\del)$ in the term
$\dv(\alpha_m(S_m^\del)\nabla S_m^\del)$
by its mean value. That is, as suggested in \cite{Arb-simpl}, we consider a constant $\overline{\alpha}_m > 0$
such that
\begin{equation}
	\label{psi_m-def}
	\overline{\alpha}_m = \int_0^1 \alpha_m(s)\, ds  \thickapprox \alpha_m(S_m^\del),
\end{equation}
and we replace the imbibition equation \eqref{imb-eqn} by its linearized version
\begin{equation}
	\label{S_m_L_delta}
	\left\{
	\begin{array}[c]{ll}
		\ds
		\Phi_m \frac{\pt \widehat{S}_{m}^{\del}}{\pt t}\, - \del^2 k_m\, \overline{\alpha}_m \Delta_y \widehat{S}_{m}^{\del} = 0
		\quad {\rm in}\,\, Y_m^{\del}, \\
		\widehat{S}_{m}^{\del}(x, y, t) = \mathcal{P}(S_f^{\del}(x,t)) \quad {\rm on}\,\, \Gamma^{\del}, \\
		\widehat{S}_{m}^{\del}(x, y, 0) = S_m^0(x) \quad {\rm in}\,\,  Y_m^{\del}.
	\end{array}
	\right.
\end{equation}
Now instead of calculating the matrix-fracture source terms defined in \eqref{H-4}
by the solution of the original imbibition equation (\ref{imb-eqn}) one can
calculate these terms using the linearized imbibition equation (\ref{S_m_L_delta}).
In this case we denote {them by}:
$\widehat{\EuScript Q}_w^{\del} = - \widehat{\EuScript Q}_n^{\del}$.
Due to {the} linearity, these terms can be expressed as a convolution. Namely,
\begin{align*}
	\widehat{\EuScript Q}_{w}^{\del}(x,t) =
	-\frac{\pt}{\pt t}\, \int_0^t {\cal K}_m^\delta(t-u)({\mathcal{P}(S_f^{\del}(x,u)) - \mathcal{P}(S_f^{0}(x))})\, du\\
	= - \int_0^t {\cal K}_m^\delta(t-u)\frac{\pt}{\pt t}\mathcal{P}(S_f^{\del}(x,u))\, du
	= - \widehat{\EuScript Q}_{n}^{\del}(x,t),
\end{align*}
where the kernel ${\cal K}_m^\delta(t)$ can easily be calculated (see, e.g., \cite{Chavant-dp}).
Due to simplicity of the domain $Y_m^\delta$ the kernel ${\cal K}_m^\delta$ can be expanded in a function 
series, truncated at some point and used in numerical calculations. In this way, it enables to avoid the resolution 
of the local problem numerically. However, the coupling between the local and global problems is still present
through  the kernel  ${\cal K}_m^\delta(t)$.
We will not give details of this approach since we are interested in the case of thin fractures.
In section~\ref{sec:decoupling} we will decouple the local and the global
problem by passing to the limit as $\delta\to 0$ and obtain a convolution form of the matrix-fracture source terms
with the  kernel given explicitly. 

We propose now a more general way to linearize
the  imbibition equation (\ref{imb-eqn}). Namely, let us consider the
following boundary value problem:
\begin{equation}
	\label{imb-eqn-lin}
	\left\{
	\begin{array}[c]{ll}
		\ds
		\Phi_m \frac{\pt \widetilde{S}_{m}^{\del}}{\pt t}\, - \del^2 k_m\, \widehat{\alpha}_m^\delta(x,t) \Delta_y \widetilde{S}_{m}^{\del} = 0
		\quad {\rm in}\,\,  Y_m^{\del}, \\
		\widetilde{S}_{m}^{\del}(x, y, t) =  \mathcal{P}(S_f^{\del}(x,t)) \quad {\rm on}\,\, \Gamma^{\del}, \\
		\widetilde{S}_{m}^{\del}(x, y, 0) = S_m^0(x) \quad {\rm in}\,\, s Y_m^{\del},
	\end{array}
	\right.
\end{equation}
where the coefficient $\widehat{\alpha}_m^\delta(x,t)$ can be chosen in different ways. One particularly 
useful choice of the coefficient $\widehat{\alpha}_m^\delta(x,t)$ is the average of the function $\alpha_m$ 
over the range of saturation given by the boundary conditions:
\begin{equation}
	\begin{split}
		\widehat{\alpha}_m^\delta(x,t) &= \int_{\mathcal{P}(S_{\rm min}^{\del}(x,t))}^{\mathcal{P}(S_{\rm max}^{\del}(x,t))} \alpha_m(s)\, ds /
		(\mathcal{P}(S_{\rm max}^{\del}(x,t)) - \mathcal{P}(S_{\rm min}^{\del}(x,t)))  \\
		&\text{where}\quad S_{\rm min}^{\del}(x,t) = \min_{0\leq s\leq t} S_f^{\del}(x,s),\quad 
		S_{\rm max}^{\del}(x,t) = \max_{0\leq s\leq t} S_f^{\del}(x,s).
	\end{split}
	\label{2nd-approx-v2-gen}
\end{equation}
Then using the function $\beta_m$ we can express $\widehat{\alpha}_m^\delta$ as follows:
\begin{align*}
	\widehat{\alpha}_m^\delta(x,t) &=\frac{ \beta_m(\mathcal{P}(S_{\rm max}^{\del}(x,t))) -\beta_m(\mathcal{P}(S_{\rm min}^{\del}(x,t))) }{
		\mathcal{P}(S_{\rm max}^{\del}(x,t))- \mathcal{P}(S_{\rm min}^{\del}(x,t))}.
\end{align*}
Note that in the case of increasing (injection of water in oil saturated media) or decreasing  fracture saturation we have 
a more simple expression for  $\widehat{\alpha}_m^\delta$:
\begin{equation}
	\begin{split}
		\widehat{\alpha}_m^\delta(x,t) &= \int_{\mathcal{P}(S_f^{\del}(x,0))}^{\mathcal{P}(S_f^{\del}(x,t))} \alpha_m(s)\, ds /
		(\mathcal{P}(S_f^{\del}(x,t)) - \mathcal{P}(S_f^{\del}(x,0)))  \\
		&=\frac{ \beta_m(\mathcal{P}(S_f^{\del}(x,t))) -\beta_m(\mathcal{P}(S_f^{\del}(x,0))) }{
			\mathcal{P}(S_f^{\del}(x,t))- \mathcal{P}(S_f^{\del}(x,0))}. 
	\end{split}
	\label{2nd-approx-v2}
\end{equation}

With any choice of the function $\widehat{\alpha}_m^\delta(x,t)\geq 0$ problem (\ref{imb-eqn-lin})
can be reduced to (\ref{S_m_L_delta}) with $\overline{\alpha}_m =1$, 
and thus to an equation with constant coefficients,
by the change of the time variable
\begin{equation}
	\label{tau-def}
	\tau_x^\delta(t) \eqdef \int\limits_0^t \widehat{\alpha}_m^\delta(x,s)\, ds
\end{equation}
and passing to a new unknown function  $\widehat{S}_m^\delta$ defined by
\begin{equation}
	\widetilde{S}_{m}^{\del}(x,y,t) = \widehat{S}_{m}^{\del}(x,y,\tau_x^\delta(t)).\label{change-of-var}
\end{equation}
One can easily check that the function $\widehat{S}_m^\delta$ is a solution to problem (\ref{S_m_L_delta}) with $\overline{\alpha}_m =1$ and 
corresponding boundary condition. Notice 
that the function $t\mapsto\tau_x^\delta(t)$ in \eqref{tau-def} is invertible except on the time intervals   where 
$\widehat{\alpha}_m^\delta(x,t) = 0$.  On these intervals, in the case of  $\widehat{\alpha}_m^\delta$ given by (\ref{2nd-approx-v2}),
the solution  $\widetilde{S}_{m}^{\del}$ of (\ref{imb-eqn-lin})
and the boundary function $S_f^{\del}$ do not depend on time,
making  (\ref{change-of-var}) consistent on these intervals and making
$\widehat{S}_m^\delta$ well defined. 

From the definition of the matrix-fracture  source terms \eqref{H-4}, using our linearized imbibition equation (\ref{imb-eqn-lin})
and the change of variables (\ref{tau-def}), (\ref{change-of-var}) we obtain  the following (linearized) form
of the matrix source term:
\begin{equation}
	\label{source-2}
	\widetilde{\EuScript Q}_w^{\del}(x,t) = - \widehat{\alpha}_m^\delta(x,t)  \frac{\Phi_m}{\mbox{vol}(Y_m^{\del})} 
	\int\limits_{Y_m^{\del}} \frac{\pt \widehat{S}_m^{\del}}{\pt \tau}(x, y, \tau_x^\delta(t))\, dy 
	= \widehat{\alpha}_m^\delta(x,t) \widehat{\EuScript Q}_w^{\del}(x,\tau_x^\delta(t)),
\end{equation}
where $ \widehat{\EuScript Q}_w^{\del}$ is the matrix-fracture  source term produced by linear imbibition equation 
(\ref{S_m_L_delta}). 

{Since} the term $\widehat{\EuScript Q}_w^{\del}(x,\tau)$ can be expressed in
{the convolution form with the kernel that is easily calculable, then}
we can {exclude the} calculation of {the} solution of the local problem
(\ref{imb-eqn-lin}) just as in the case of {the} constant approximation (\ref{psi_m-def}),
(\ref{S_m_L_delta}).

{\bf Remark}. {\em Another possible linearization consists in rewriting equation \eqref{imb-eqn} in terms of 
	the function $U^\delta =  \beta_m(S_m^{\del})$ with nonlinearity now appearing in the time derivative term.
	One can then apply the same type of linearizations as above which lead 
	to different fracture-matrix source terms. In the numerical simulations with this kind of linearization
	we have observed different results but without net improvement over linearization we have presented above,
	so we do not consider that kind of linearization further. }

\section{Simplified dual porosity model}
\label{sec:decoupling}

In dual porosity model given by (\ref{H-1}), (\ref{H-4}) and (\ref{imb-eqn}) there remains a 
small parameter $\delta$ which measures relative thickness of the fractures. By letting this small
parameter to zero one can obtain full decoupling of the local and the global problem. 
In the case of the two-phase flow that problem was studied in \cite{JPV}; see the references therein 
for the case of one-phase flow.

It is shown in \cite{JPV} that in the limit  $\delta\to 0$ the effective fracture equations (\ref{H-1})
coupled with the linearized imbibition equation (\ref{S_m_L_delta}) reduces to   the system:
\begin{equation}
	\label{H-1a}
	\left\{
	\begin{array}[c]{ll}
		\ds
		\Phi_f   \frac{\pt S}{\pt t} -
		{\rm div}\, \bigg( k^*\lm_{w,f} (S) \gr P_{w} \bigg) = {\EuScript Q}_w \\
		\ds
		- \Phi_f \frac{\pt S}{\pt t} -
		{\rm div}\, \bigg(k^* \lm_{n,f} (S) \gr P_{n} \bigg) = {\EuScript Q}_n, \\
		P_{c,f}(S) = P_{n} - P_{w},
	\end{array}
	\right.
\end{equation}
where $S$, $P_n$ and $P_w$ are effective variables in the system of fractures, $k^* = k_f(d-1)/d$ and source terms ${\EuScript Q}_w = -{\EuScript Q}_n$ 
are given by 
\begin{align}
	{\EuScript Q}_w(x,t) = -  C_m \frac{\pt}{\pt t}\, \int_0^t \frac{\mathcal{P}(S(x,u)) - \mathcal{P}(S(x,0))}{\sqrt{t-u}}\, du, 
	\label{Q-modelI}\\
	C_m = 2d\sqrt{{\Phi_m k_m\overline{\alpha}_m}/{\pi}}.
	\label{Q-modelI-C}
\end{align}
We note that in  this final reduction  the effective porosity is replaced by the fracture porosity,
the effective permeability is reduced fracture permeability, and 
the matrix-fracture source  term ${\EuScript Q}_w^{\del}(t)$ reduces to the convolution expression 
(\ref{Q-modelI}) with the kernel ${\cal K}(t) = C_m/\sqrt{t}$. In the model (\ref{H-1a}), (\ref{Q-modelI}) 
effective permeability and the matrix-fracture source  term are given explicitly and there is no need
to solve local problems. The local and the global problems are completely decoupled. 
We should also mention that in the system  (\ref{H-1}), (\ref{S_m_L_delta}) the porosity, 
permeability and the  matrix-fracture source  term are proportional to $\delta$, and equations (\ref{H-1a}), (\ref{Q-modelI})
are obtained after division by  $\delta$. For the source terms we therefore have 
$\widehat{\EuScript Q}_{w}^{\del} = \delta {\EuScript Q}_w + o(\delta)$ (see \cite{JPV} for more details).

In the case of variable linearization of the imbibition equation given by (\ref{imb-eqn-lin})
we can use the matrix-fracture source representation (\ref{source-2}) to compute, at least formally, 
the asymptotic behavior of $\widetilde{\EuScript Q}_{w}^{\del}$ when $\delta\to 0$. To this end we 
assume that the boundary function $S_f^\delta(x,t)$ converges to some $S_f(x,t)$ leading to convergence 
$ \widehat{\alpha}_m^\delta(x,t)\to  \widehat{\alpha}_m(x,t)$ and $\tau_x^{\del}(t)\to \tau_x(t)$, where 
obviously
\begin{equation}
	\begin{split}
		\widehat{\alpha}_m(x,t)
		&=\frac{ \beta_m(\mathcal{P}(S_f(x,t))) -\beta_m(\mathcal{P}(S_f(x,0))) }{
			\mathcal{P}(S_f(x,t))- \mathcal{P}(S_f(x,0))},
	\end{split}
	\label{2nd-approx-v2-limit}
\end{equation}
and
\begin{equation}
	\label{tau-def-limit}
	\tau_x(t) \eqdef \int\limits_0^t \widehat{\alpha}_m(x,s)\, ds.
\end{equation}
Now formula (\ref{source-2}) gives
\begin{equation*}
	\widetilde{\EuScript Q}_w^{\del}(x,t)= \widehat{\alpha}_m^\delta(x,t) \widehat{\EuScript Q}_w^{\del}(x,\tau_x^\delta(t)),
\end{equation*}
where by asymptotic expansion of the matrix-fracture source term in the case of constant linearization we have 
\begin{equation*}
	\widehat{\EuScript Q}_{w}^{\del}(x,\tau) =
	-\del C_m \frac{\pt}{\pt \tau}\, \int_0^{\tau} 
	\frac{\mathcal{P}(S_f(x,(\tau_x)^{-1}(u))) - \mathcal{P}(S_f(x,0)) }{\sqrt{\tau -u}}\, du +o(\delta),
\end{equation*}
and $C_m = 2d\sqrt{{\Phi_m k_m }/{\pi}}$.
Therefore we get   
\begin{equation}
	\widetilde{\EuScript Q}_w^{\del}(x,t) 
	= - \del C_m\frac{\pt}{\pt t}\, \int_0^{\tau_x(t)}
	\frac{\mathcal{P}(S_f(x,(\tau_x)^{-1}(u))) - \mathcal{P}(S_f(x,0)) }{\sqrt{\tau_x(t) -u}}\, du+o(\delta).
	\label{Q-develop-tmp}
\end{equation}
Further change of variables  $u=\tau_x(s)$ gives
\begin{equation*}
	\widetilde{\EuScript Q}_w^{\del}(x,t) 
	= - \del C_m\frac{\pt}{\pt t}\, \int_0^{t}\frac{\mathcal{P}(S_f(x,s)) - \mathcal{P}(S_f(x,0)) }
	{\sqrt{\tau_x(t) -\tau_x(s)}}\widehat{\alpha}_m(x,s)\, ds +o(\delta).
\end{equation*}

We can now conclude that in the limit  $\delta\to 0$ the effective fracture equations (\ref{H-1})
coupled with the linearized imbibition equation (\ref{imb-eqn-lin}) reduces to   the system  (\ref{H-1a})
with the effective matrix-fracture source term  ${\EuScript Q}_w =-{\EuScript Q}_n$ given by
\begin{equation}
	{\EuScript Q}_w(x,t) 
	= - C_m\frac{\pt}{\pt t}\, \int_0^{t}\frac{\beta_m(\mathcal{P}(S_f(x,s))) -\beta_m(\mathcal{P}(S_f(x,0))}
	{\sqrt{\tau_x(t) -\tau_x(s)}}\, ds,
	\label{Q-develop-tmp-2}
\end{equation}
where $\tau_x$ is given by (\ref{2nd-approx-v2-limit}), (\ref{tau-def-limit}) and  $C_m = 2d\sqrt{{\Phi_m k_m }/{\pi}}$.
We have obtained ${Q}_w$ in (\ref{Q-develop-tmp-2}) by leaving out small $o(\delta)$ terms and dividing by $\delta$ in  (\ref{Q-develop-tmp}).

\section{Numerical comparison of  matrix-fracture exchange terms}
\label{sec:num_comp}

In this section we will compare by means of numerical simulation the matrix-fracture exchange term calculated by solving
nonlinear imbibition equation (\ref{imb-eqn}) to the matrix-fracture exchange terms given by different linearisation procedures. 
Since we do not consider the whole two-phase flow simulation in this section we will impose an artificial boundary condition on 
the block which corresponds to water injection into oil saturated media.

In this example we consider porous block $Y = (0,1-\delta)^{d}$  and we will consider relative fracture thickness 
$\delta=0.3, 0.1, 0.01, 0.001$.  
For fracture permeability we will take fixed value $k_f = 10^{-13}$ m$^2$ and we take $k_m = k_f$. 
This means that the matrix permeability in the scaled matrix block $Y$ is equal to $\delta^2 k_f$. 
Furthermore, we take  the matrix porosity $\Phi_m = 0.35$ and 
fluid viscosities $\mu_w = 10^{-3}$ Pa$\cdot$s and $\mu_n = 2\cdot 10^{-3}$ Pa$\cdot$s.

Van Genuchten-Mualem model is used to express the capillary pressure functions 
and the relative permeabilities. We have 
\begin{align*}
	P_c(S_{w})     &= P_r (S_{w}^{-1/m} -1)^{1/n},\\
	k_{rw}(S_{w})  = \sqrt{S_{w}} [ 1 - (1 - S_{w}^{1/m})^m]^2,&\quad
	k_{rg}(S_{w})  = \sqrt{1-S_{w}} (1 - S_{w}^{1/m})^{2m},
\end{align*}
where  $S_{w}$ is effective wetting phase saturation and the  residual saturations are taken to be zero;
$P_r$ is reference pressure for Van Genuchten law  and
$n > 1$, $m\geq 0$ are such that $m = 1-1/n$.

\begin{figure}[h]
	\centering
	\input{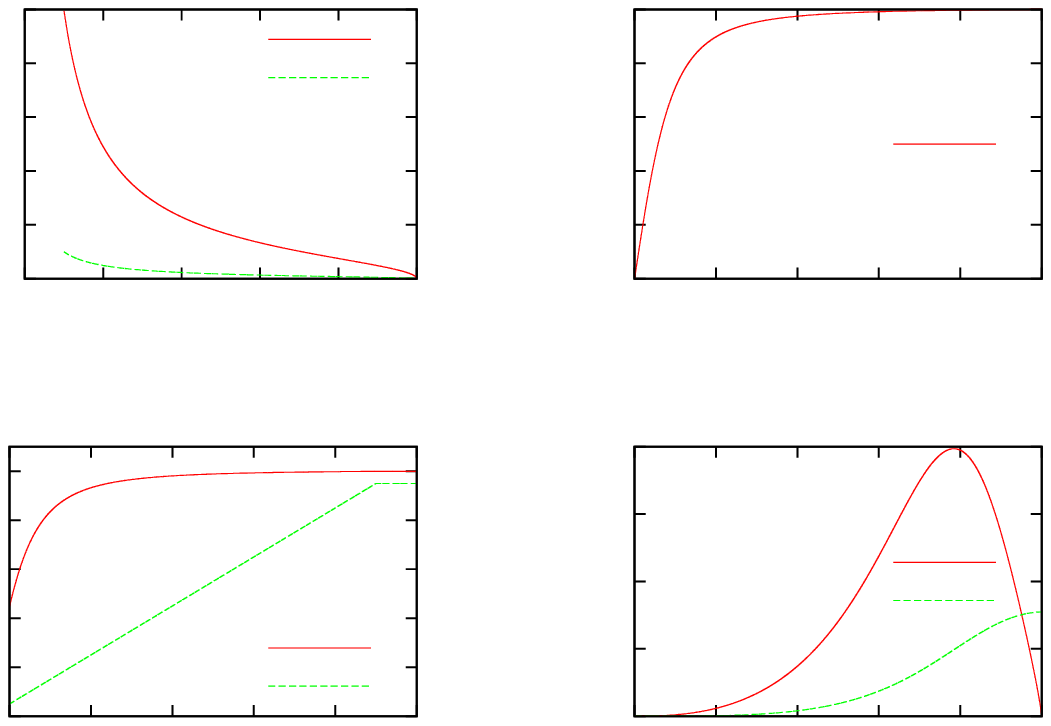}
	\caption{Functions used in Simulation 1. Van Genuchten parameters are: $P_r = 1$ bar and $n=2$ 
		in the matrix and $P_r = 0.1$ bar and $n=2$ in the fractures.}
	\label{Fig:f1}
\end{figure}

Evolution of the fracture saturation is given by 
\begin{align}
	S_w^f(t) =  0.05 +\min(t/10,0.9),\label{fracture-sat}
\end{align}
and is shown on Figure~\ref{Fig:f1} together with the boundary condition ${\cal P}(S_w^f(t))$. 
On the same figure we also show the capillary pressure functions in the fracture $P_{c,f}$ and
in the matrix $P_{c,m}$ as well as the functions $\alpha_m$, $\beta_m$ and ${\cal P}$.

For numerical resolution of the imbibition equation we need  a grid that is well 
adapted for resolving the boundary layers that governs the solution. For that purpose we have 
adopted meshes of Bakhvalov type (see \cite{Linss}). The adequacy of chosen parameters of 
the Bachvalov  grid is verified in two ways. First, in 1-D we  compared numeric solution
to the analytic solution which is easy to
calculate in the linear cases. Secondly, since the mass transfer term can be calculated 
by volume integration and by boundary integration we refined the grid up to point where 
the two methods give results that differ no more than 1 \%.

\begin{figure}[t]
	\centering
    \input{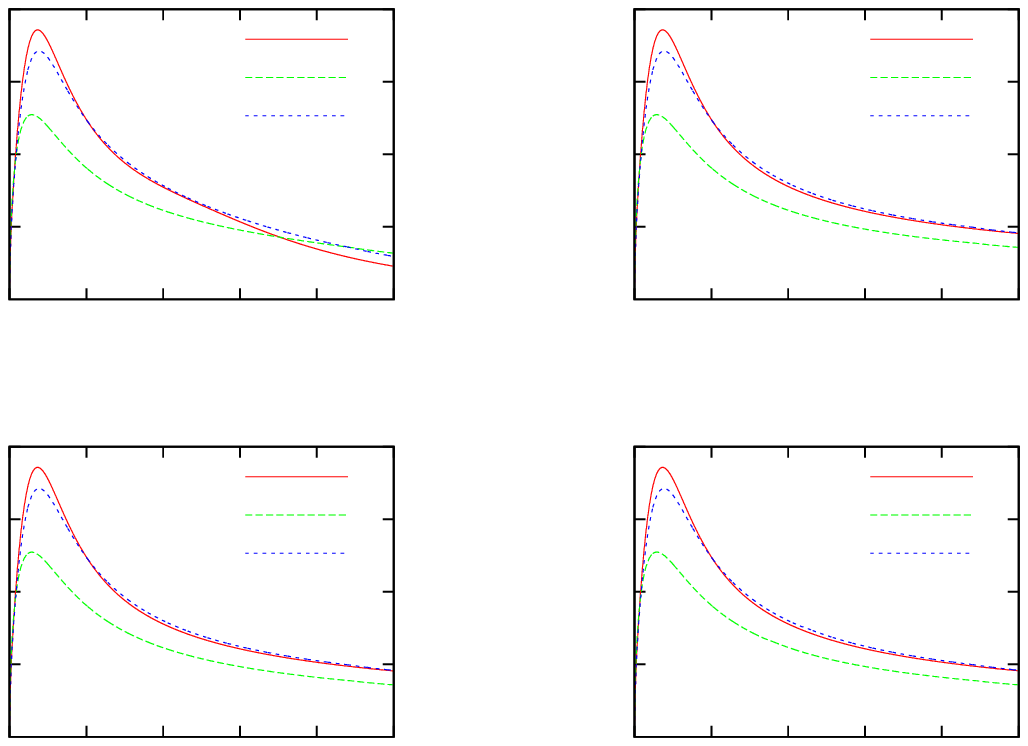}
	\caption{Evolution of the matrix-fracture exchange term  ${\EuScript Q}_w^{\del}(t)$ diveded by 
		$\delta$ for different values of $\delta$, for the nonlinear and two linear models. On $x$-axis is given time in days. }
	\label{Fig:f6}
\end{figure}

	On Figure~\ref{Fig:f6} we present time evolution of the matrix-fracture transfer term ${\EuScript Q}_w^{\del}(t)$ 
	for the nonlinear model (\ref{imb-eqn}) (denoted by {\sf nlin})
	and for two linear models: one given by (\ref{S_m_L_delta}) (denoted by {\sf clin})
	and model (\ref{imb-eqn-lin}) with (\ref{2nd-approx-v2}) (denoted by {\sf vlin}).
Figure~\ref{Fig:f6} shows the matrix-fracture transfer term ${\EuScript Q}_w^{\del}(t)$ divided by $\delta$
for different values of 
fracture thickness $\delta$. It is shown that the approximation to nonlinear  matrix-fracture transfer term
given by  linear  model
(\ref{imb-eqn-lin}) and (\ref{2nd-approx-v2}) is much better than approximation given by model (\ref{S_m_L_delta})
and that the quality of the approximation is rather independent of the fracture size $\delta$. In fact
for small $\delta$, expression ${\EuScript Q}_w^{\del}(t)/\delta$ becomes quickly actually independent of $\delta$,
which confirms theoretical result in \cite{JPV}.

\begin{figure}[p]
	\centering
    \input{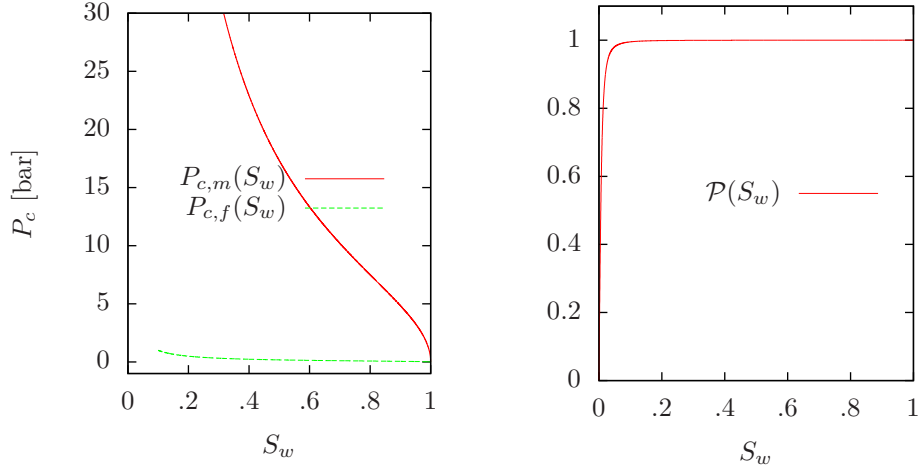}
	\caption{The capillary pressure functions and saturation transfer 
		curve in the case of large difference in the matrix and the fracture. 
		Van Genuchten parameters are: $P_r = 10$ bar and $n=2$ 
		in the matrix and $P_r = 0.1$ bar and $n=2$ in the fractures.}
	\label{Fig:f7}
\end{figure}

The matrix-fracture exchange term is strongly influenced by difference between the capillary pressure curves in 
the matrix and the fracture. If the difference between the two curves is large the saturation transfer 
function ${\cal P}$ will have strong derivative near $S_w=0$ and it will be almost constant in the rest of the domain. 
This will strongly influence the boundary condition for the imbibition equation. In Figure~\ref{Fig:f7}
we show the case of the van Genuchten capillary pressure functions with the parameters $n=2$ and $P_r=10$ bars 
in the matrix and with  $n=2$ and $P_r=0.1$ bars in the fracture. In the other extreme, where the two 
capillary pressure functions are mutually equal, the saturation transfer 
function ${\cal P}$ is linear. 

\begin{figure}[p]
	\centering
    \input{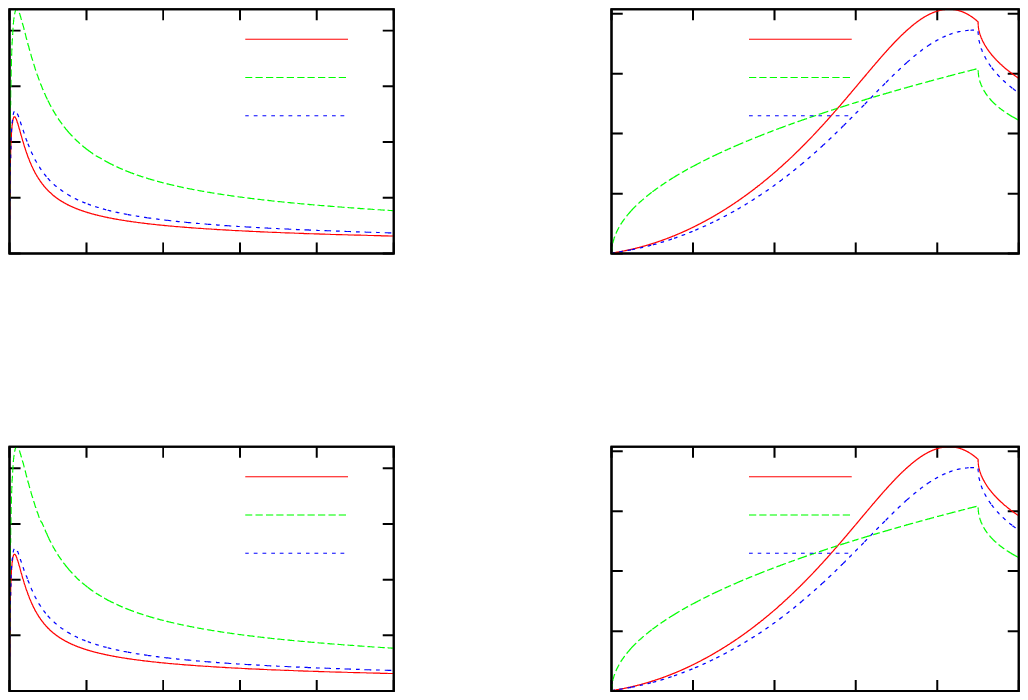}
	\caption{In the left column strong difference in the matrix and fracture $P_c$ functions. 
		Parameters are given in Figure~\ref{Fig:f7};
		in the right column case of equal $P_c$ functions. In the first row $\delta =10^{-2}$ and in the second row
		$\delta =10^{-3}$.}
	\label{Fig:f8}
\end{figure}

In Figure~\ref{Fig:f8} we compare the matrix-fracture exchange term ${\EuScript Q}_w^{\del}(t)$
in the case of strongly different capillary pressure
functions shown in Figure~\ref{Fig:f7}, and the case of equal capillary pressure functions ($P_r = 1$ bar, $n=2$) 
in the matrix and the fracture;
in both case $\del = 10^{-3}$. We see in both cases that better approximation to the nonlinear  matrix-fracture exchange term
is again given by {\sf vlin} curve, that is by the  model given by (\ref{imb-eqn-lin}) and (\ref{2nd-approx-v2}).
Simpler {\sf clin} approximation given by model (\ref{S_m_L_delta}) gives in all cases less good approximation.

\begin{figure}[p]
	\centering
	\input{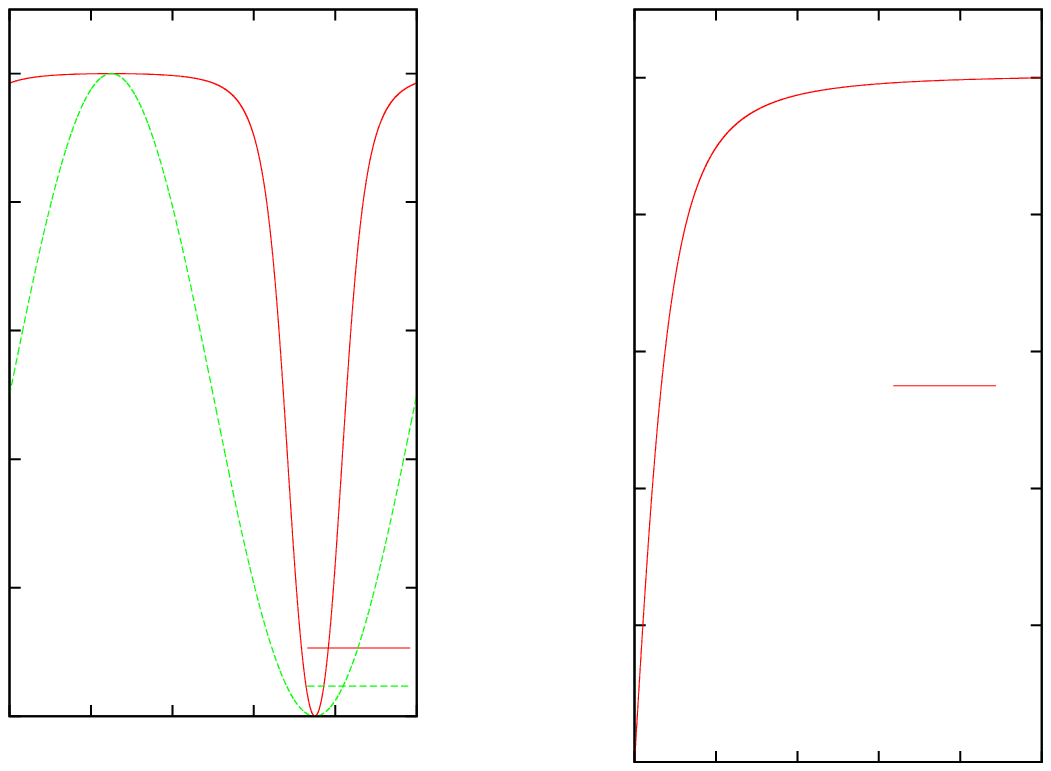}
	\caption{Non monotone boundary saturation. }
	\label{Fig:ff5}
\end{figure}

\begin{figure}[p]
	\centering
	\input{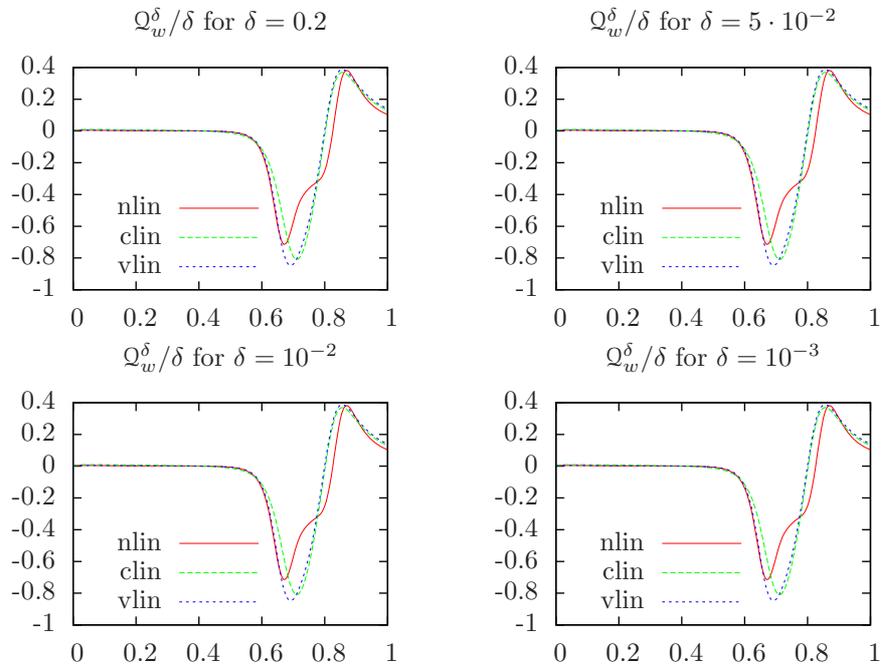}
	\caption{Matrix-fracture exchange term in the case of non monotone boundary condition and  $\delta =   10^{-3}$. }
	\label{Fig:ff6}
\end{figure}

Finally we chose an example in which the fracture saturation is not monotone. In that case we need to use 
the definition of the coefficient $\widehat{\alpha}_m^\delta(x,t)$ given by (\ref{2nd-approx-v2-gen}).
Keeping all the other parameters as before we change only the boundary conditions which is now given by 
\begin{equation}
	S_w^f(t) = 0.5 + 0.5 \sin(\pi t/5),
	\label{fracture-sat-nm}
\end{equation}
and the simulation time is, as before, 10 days. 
The function (\ref{fracture-sat-nm}) and corresponding boundary condition ${\cap P}(S_w^f(t))$ are shown on 
Figure~\ref{Fig:ff5}. The  matrix-fracture exchange term ${\EuScript Q}_w^{\del}(t)$ is shown on  Figure~\ref{Fig:ff6}.
We see that after lost of the monotonicity of the boundary condition the  matrix-fracture exchange term given by
the variable linearization 
looses its precision but stays comparable to the constant linearization version of the matrix-fracture exchange term.

\section{Discretization of dual porosity model I}
\label{sec:discrete1}

The model (\ref{H-1a})--(\ref{Q-modelI}) will be discretized by the cell centered finite volume method on a structured grid
with the two-point flux approximation. First we present the time discretization. 


Assume that we have a sequence of time steps: $0=t^0 < t^1 < \cdots < t^n < \cdots$ and denote   $\delta t^n = t^{n+1}-t^n$
and also $I_n = (t_{n-1},t_{n}]$. All unknowns are supposed  to be 
piecewise constant in time, such that  $ S(x,t) = \sum_k S^{k}(x) \chi_{I_k}(t)$, where 
$S^{k}(x) = S(x,t_k)$, and similarly for other variables. Implicit Euler discretization gives for $t\in I_{n+1}$,
\begin{align*}
	\Phi_f \frac{S^{n+1} - S^{n}}{\delta t^n}-{\rm div}\, \bigg(  \lm_{w,f} (S^{n+1})k^* \gr P_{w}^{n+1} \bigg) = Q^{n+1/2} , \\
	- \Phi_f \frac{S^{n+1} - S^{n}}{\delta t^n} - {\rm div}\, \bigg( \lm_{n,f} (S^{n+1}) k^* \gr P_{n}^{n+1} \bigg) = -Q^{n+1/2}. 
\end{align*}
The source term is discretized in the following way:
\begin{align*}
	{Q}^{n+1/2}  =       -\frac{C_m}{\delta t^n} 
	\Big(&\int_0^{t_{n+1}} \sum_{k=1}^{n+1} \frac{\mathcal{P}(S^{k}) - \mathcal{P}(S^{0})}
	{\sqrt{t_{n+1}-s}}\chi_{I_k}(s)\, ds\\
	&- \int_0^{t_{n}}    \sum_{k=1}^{n} \frac{\mathcal{P}(S^{k}) - \mathcal{P}(S^{0})}
	{\sqrt{t_{n}-s}}\chi_{I_k}(s)\, ds\Big)\\
	= - \frac{1}{\delta t^n} 
	\Big(&  \sum_{k=1}^{n+1}  (\mathcal{P}(S^{k}) - \mathcal{P}(S^{0}) )\int_{t_{k-1}}^{t_{k}}
	\frac{C_m ds}  {\sqrt{t_{n+1}-s}}\\
	&-    \sum_{k=1}^{n}  ( \mathcal{P}(S^{k}) - \mathcal{P}(S^{0})  )  \int_{t_{k-1}}^{t_{k}}   \frac{C_m ds}
	{\sqrt{t_{n}-s}}\Big)\\
	= - \frac{C_m}{\delta t^n} 
	\Big(&  \sum_{k=1}^{n+1}  (\mathcal{P}(S^{k}) - \mathcal{P}(S^{0}) ) I^{n+1}_k 
	-    \sum_{k=1}^{n}  ( \mathcal{P}(S^{k}) - \mathcal{P}(S^{0})  )  I^{n}_{k}\Big)\\
	= - \frac{C_m}{\delta t^n} 
	\Big(&  (\mathcal{P}(S^{n+1}) - \mathcal{P}(S^{0}) ) I^{n+1}_{n+1}
	+ \sum_{k=1}^{n}  (\mathcal{P}(S^{k}) - \mathcal{P}(S^{0}) ) (I^{n+1}_{k} -  I^{n}_{k})\Big)
\end{align*}
where we  denoted  for $1\leq k\leq n$,
\begin{align*}
	I^{n}_{k} =  \int_{t_{k-1}}^{t_{k}}   \frac{C_m ds} {\sqrt{t_{n}-s}} = 2C_m(\sqrt{t_n - t_{k-1}} - \sqrt{t_n - t_k})
	= \frac{\delta t^{k-1}}{\sqrt{t_n - t_{k-1}}+ \sqrt{t_n - t_k}}.
\end{align*}
Obviously, $I^{n+1}_{n+1} = 2C_m \sqrt{\delta t^n}$. If the time grid is equidistant, then we have $I^{n+1}_{k+1}=I^{n}_{k} = J_{n-k}$,
since
\begin{align*}
	\int_{t_{k}}^{t_{k+1}}   \frac{ds}  {\sqrt{t_{n+1}-s}} =  \int_{t_{k-1}}^{t_{k}}   \frac{ds} {\sqrt{t_{n}-s}}
	= \frac{2\sqrt{\delta t}}{\sqrt{n-k+1}+\sqrt{n-k}},\quad
	J_l = \frac{2C_m\sqrt{\delta t}}{\sqrt{l+1}+\sqrt{l}},
\end{align*}
leading to a convolution-like representation:
\begin{align*}
	{Q}^{n+1/2}  = -\frac{1}{\delta t} \sum_{k=0}^{n}  [\mathcal{P}(S^{k+1}) - \mathcal{P}(S^{k}) ] J_{n-k}. 
\end{align*}

Generally we have:
\begin{equation}
	\begin{split}
		\Phi_f \frac{S^{n+1}}{\delta t^n} &+ \frac{2C_m}{\sqrt{\delta t^n}} \mathcal{P}(S^{n+1})
		-{\rm div}\, \bigg( \lm_{w,f} (S^{n+1}) k^*\gr P_{w}^{n+1} \bigg)\\
		&=  \Phi_f \frac{S^{n}}{\delta t^n} + \frac{1}{\delta t^n}{\cal F}^n 
	\end{split}
	\label{disc-w-tipI}
\end{equation}
\begin{equation}
	\begin{split}
		- \Phi_f \frac{S^{n+1}}{\delta t^n} 
		&- \frac{2C_m}{\sqrt{\delta t^n}} \mathcal{P}(S^{n+1}) 
		- {\rm div}\, \bigg(\lm_{n,f} (S^{n+1}) k^*\gr P_{n}^{n+1} \bigg)\\
		& = - \Phi_f \frac{ S^{n}}{\delta t^n}  - \frac{1}{\delta t^n}{\cal F}^n
	\end{split}
	\label{disc-n-tipI}
\end{equation}
where, for $n>0$,
\begin{align}
	\label{Fn-tipI}
	{\cal F}^n & =   
	\mathcal{P}(S^{0})  I^{n+1}_{n+1}
	- \sum_{k=1}^{n}  (\mathcal{P}(S^{k}) - \mathcal{P}(S^{0}) ) (I^{n+1}_{k} -  I^{n}_{k})\\
	& = \mathcal{P}(S^{0}) ( I^{n+1}_{n+1}+\sum_{k=1}^{n} (I^{n+1}_{k} -  I^{n}_{k}) ) 
	- \sum_{k=1}^{n}  \mathcal{P}(S^{k}) (I^{n+1}_{k} -  I^{n}_{k}).
\end{align}
Also, note that in the case $n=0$ we have,
\begin{align*}
	{Q}_w^{1/2} &=  -\frac{C_m}{\delta t^0} 
	\Big(\int_0^{t_{1}}  \frac{\mathcal{P}(S^{1}) - \mathcal{P}(S^{0})} {\sqrt{t_{1}-s}}\, ds   - 0\Big)\\
	&= - \frac{C_m}{\delta t^0} 
	(\mathcal{P}(S^{1}) - \mathcal{P}(S^{0}) )\int_{0}^{t_{1}}\frac{ds}  {\sqrt{t_{1}-s}}.
\end{align*}
Therefore, for $n=0$ we have 
\begin{align*}
	{\cal F}^0 =  \mathcal{P}(S^{0})  I^{1}_{1} = 2 \sqrt{\delta t^0}  \mathcal{P}(S^{0}),
\end{align*}
and (\ref{Fn-tipI}) holds also for $n=0$.

Let us denote $ D_k^n =  I^{n}_{k} -I^{n+1}_{k}$ for $k=1,\ldots, n$. We have 
\begin{align*}
	D_k^n =  I^{n}_{k} -I^{n+1}_{k} = \frac{C_m \delta t^{k-1}}{\sqrt{t_n - t_{k-1}}+ \sqrt{t_n - t_k}} - 
	\frac{C_m \delta t^{k-1}}{\sqrt{t_{n+1} - t_{k-1}}+ \sqrt{t_{n+1} - t_k}} > 0
\end{align*}
since the function
\begin{align*}
	\omega(t) = \frac{C_m}{\sqrt{t}}
\end{align*}
is monotone decreasing. We also introduce 
\begin{align*}
	D^n_0 =  I^{n+1}_{n+1} + \sum_{k=1}^{n} (I^{n+1}_{k} -  I^{n}_{k}) =
	I_1^{n+1} +\sum_{k=1}^{n} (I^{n+1}_{k+1} - I^n_k) .
\end{align*}
Let us note that in the equidistant time stepping we have $I^{n+1}_{k+1} - I^n_k=0$ and then $D^n_0 = I_1^{n+1} >0$. 
In the non equidistant case the terms $I^{n+1}_{k+1} - I^n_k$ can have any sign, so we will introduce the assumption 
that the time discretization is such that 
\begin{align}
	D^n_0 >0.
	\label{time-disc-cond-0}
\end{align}
This will always be the case if the time stepping is close to the equidistant one.

With introduced notation we can write 
\begin{align}
	\label{Fn-tipI-new}
	{\cal F}^n & = \sum_{k=0}^{n}  D_k^n \mathcal{P}(S^{k})  
\end{align}
with $D_k^n >0$ for $k=0,1,\ldots,n$.

Let us also note that 
\begin{align}
	\sum_{k=0}^{n}  D_k^n =\sum_{k=1}^{n} (I^{n}_{k} -I^{n+1}_{k}) + I^{n+1}_{n+1} + \sum_{k=1}^{n} (I^{n+1}_{k} -  I^{n}_{k}) 
	= I^{n+1}_{n+1}= 2C_m \sqrt{\delta t^n}.
	\label{sum-dnk}
\end{align}


We use standard finite volume discretization of the two phase system written in the phase formulation 
\cite{EHM-petroleum} (see \cite{EGH-book} for notations): 

\begin{equation}
	\begin{split}
		\Phi_{f,K} \frac{S^{n+1}_K}{\delta t^n} &+ \frac{2C_m}{\sqrt{\delta t^n}} \mathcal{P}(S^{n+1}_K)
		- \sum_{L\in {\cal N}_K}\tau_{K|L}  k^*_{K|L}  \lm_{w,f,K|L}^{n+1} \delta_{K,L}^{n+1}(P_{w})\\
		&=  \Phi_f \frac{S^{n}_K}{\delta t^n} + \frac{1}{\delta t^n}{\cal F}^n_K 
	\end{split}
	\label{disc-w-tipI-full}
\end{equation}
\begin{equation}
	\begin{split}
		- \Phi_{f,K} \frac{S^{n+1}_K}{\delta t^n} 
		&- \frac{2C_m}{\sqrt{\delta t^n}} \mathcal{P}(S^{n+1}_K) 
		-  \sum_{L\in {\cal N}_K}\tau_{K|L}  k^*_{K|L}  \lm_{n,f,K|L}^{n+1} \delta_{K,L}^{n+1}(P_{n})\\
		& = - \Phi_{f,K} \frac{S^{n}_K}{\delta t^n}  - \frac{1}{\delta t^n} {\cal F}^n_K
	\end{split}
	\label{disc-n-tipI-full}
\end{equation}

In this discretization we use phase by phase upstream choice: the value of the mobility of each phase on the edge $K|L$ is
determined by the sign of the  difference of the discrete phase pressure. 

\begin{align}
	\lm_{w,f,K|L}^{n+1} =\lm_{w,f}( S_{w,K|L}^{n+1}),\quad  \lm_{n,f,K|L}^{n+1}=\lm_{n,f}( S_{n,K|L}^{n+1}),
	\label{upwind:1}
\end{align}
with 
\begin{equation}
	\begin{split}
		&   S_{w,K|L}^{n+1} = \begin{cases} S_{K}^{n+1} & \text{if } (K,L)\in {\cal E}_w^{n+1}\\  S_{L}^{n+1} & \text{otherwise},\end{cases}\quad
		S_{n,K|L}^{n+1} = \begin{cases} S_{K}^{n+1} & \text{if } (K,L)\in {\cal E}_n^{n+1}\\  S_{L}^{n+1} & \text{otherwise},\end{cases}\\
		&  \text{where } {\cal E}_w^{n+1}\text{ and } {\cal E}_n^{n+1} \text{ are two subsets of } {\cal E} \text{ such that }\\
		& \{(K,L)\in {\cal E}\colon  \delta_{K,L}^{n+1}(P_{w}) < 0\} \subset {\cal E}_w^{n+1} \subset   
		\{(K,L)\in {\cal E}\colon  \delta_{K,L}^{n+1}(P_{w}) \leq 0\} \\
		&    \{(K,L)\in {\cal E}\colon  \delta_{K,L}^{n+1}(P_{n}) < 0\} \subset {\cal E}_n^{n+1} \subset   
		\{(K,L)\in {\cal E}\colon  \delta_{K,L}^{n+1}(P_{n}) \leq 0\}
	\end{split}
	\label{upwind:2}
\end{equation}

\section{Discretization of dual porosity model II}
\label{sec:discrete2}

Second model differs from the first one only in the matrix-fracture exchange term 
which takes the form: 
\begin{equation*}
	\widetilde{\EuScript Q}_w(x,t) 
	= -  C_m\frac{\pt}{\pt t}\, \int_0^{\tau_x(t)}
	\frac{\mathcal{P}(S(x,(\tau_x)^{-1}(u))) - \mathcal{P}(S(x,0)) }{\sqrt{\tau_x(t) -u}}\, du,
\end{equation*}
where $C_m = 2\sqrt{{\Phi_m k_m}/{\pi}}$ and  $\tau_x$ is given by
\begin{equation*}
	\tau_x(t) = \int\limits_0^t \widehat{\alpha}_m(x,s)
	\, ds.
\end{equation*}
We have chosen expression for the matrix-fracture exchange term given by (\ref{Q-develop-tmp}) but it is also possible
to use other forms, for example (\ref{Q-develop-tmp-2}). 

We will discretize this model using the same approach as in the constant linearization model. 
Assume that we have a sequence of time instances $0=t^0 < t^1 <\cdots <t^n < \cdots $ and denote by $\tau_x^n = \tau_x(t^n)$.
If the saturation  $S$  is constant in time on each interval $(t^k,t^{k+1})$ then 
$\mathcal{P}(S(x,(\tau_x)^{-1}(u))) - \mathcal{P}(S(x,0)) $ is constant on each interval $(\tau_x^k,\tau_x^{k+1})$ 
and we can write
\begin{align*}
	\widetilde{\EuScript Q}_w^{n+1/2} &\approx - \frac{C_m}{\delta t^n} 
	\Big(\int_0^{\tau_x^{n+1}} \frac{\mathcal{P}(S(x,(\tau_x)^{-1}(u))) -\mathcal{P}(S(x,0))} {\sqrt{\tau_x^{n+1} - u }} \, du\\
	&\qquad\quad - \int_0^{\tau_x^{n}}\frac{\mathcal{P}(S(x,(\tau_x)^{-1}(u))) -\mathcal{P}(S(x,0))} {\sqrt{\tau_x^n - u }}  \, du\Big)\\
	& = - \frac{1}{\delta t^n}
	\Big(\sum_{k=1}^{n+1} \int_{\tau_x^{k-1}}^{\tau_x^{k}} \frac{C_m du} {\sqrt{\tau_x^{n+1} - u}} 
	(\mathcal{P}(S^k(x)) -\mathcal{P}(S^0(x)) )\\
	&\qquad\quad - \sum_{k=1}^{n} \int_{\tau_x^{k-1}}^{\tau_x^{k}} \frac{C_m du} {\sqrt{\tau_x^n - u}} 
	(\mathcal{P}(S^k(x)) -\mathcal{P}(S^0(x)) )
	\Big).
\end{align*}
As before we have 
\begin{align*}
	I^{n}_{k} = \int_{\tau_x^{k-1}}^{\tau_x^{k}} \frac{C_m du} {\sqrt{\tau_x^{n} - u}} 
	= 2C_m(\sqrt{\tau_x^{n} - \tau_x^{k-1}} - \sqrt{\tau_x^{n} -  \tau_x^{k}})
	= \frac{ 2C_m( \tau_x^{k}-  \tau_x^{k-1}) }{\sqrt{\tau_x^{n} -  \tau_x^{k-1}}+ \sqrt{\tau_x^{n} -  \tau_x^{k}}}.
\end{align*}

Note that for $s\in (t^{k-1},t^k)$ we have 
\begin{align*}
	\tau(s) = \int_0^s \widehat{\alpha}_m(u)\, du = \sum_{l=1}^{k-1} \widehat\alpha_m^l \delta t^{l-1} + \widehat\alpha_m^k (s - t^{k-1}),
\end{align*}
where $\widehat\alpha_m^l= \widehat\alpha_m(t^l)$, so that 
\begin{align*}
	\tau(t^n) - \tau(t^{k-1}) = \sum_{l=k}^{n} \widehat\alpha_m^l \delta t^{l-1} ,\quad
	\tau(t^n) - \tau(t^k) =\sum_{l=k+1}^{n} \widehat\alpha_m^l \delta t^{l-1}, \quad \tau(t^{k+1}) - \tau(t^{k}) =  \widehat\alpha_m^k \delta t^{k-1}.
\end{align*}
Therefore, we have for $k\leq n$,
\begin{align*}
	I^{n}_{k} = \frac{ 2C_m\widehat\alpha_m^k \delta t^{k-1} }{\sqrt{\sum_{l=k}^{n} \widehat\alpha_m^l \delta t^{l-1}}
		+ \sqrt{\sum_{l=k+1}^{n} \widehat\alpha_m^l \delta t^{l-1}}}.
\end{align*}
For notational simplicity we will introduce for $k\leq n$ 
\begin{align}
	U_k^n = \sum_{l=k}^{n} \widehat\alpha_m^l \delta t^{l-1}, \label{U-k-n}
\end{align}
and $U_k^n = 0$ for $k>n$. Then we can write:
\begin{align*}
	\widetilde{\EuScript Q}_w^{n+1/2} &= - \frac{2C_m}{\delta t^n}
	\Big(\sum_{k=1}^{n+1} \frac{ \widehat\alpha_m^k (\mathcal{P}(S^k) -\mathcal{P}(S^0) ) }{ \sqrt{U_{k}^{n+1}} + \sqrt{U_{k+1}^{n+1}} }
	\delta t^{k-1} 
	- \sum_{k=1}^{n}   \frac{ \widehat\alpha_m^k  (\mathcal{P}(S^k) -\mathcal{P}(S^0) ) }{\sqrt{U_{k}^{n}} + \sqrt{U_{k+1}^{n} }}
	\delta t^{k-1} 
	\Big).
\end{align*}

We finally obtain the following scheme:
\begin{equation}
	\begin{split}
		\Phi_f \frac{S^{n+1}}{\delta t^n} & + \frac{2C_m}{\delta t^n} 
		\sum_{k=1}^{n+1}  \frac{\widehat\alpha_m^k (\mathcal{P}(S^k) -\mathcal{P}(S^0) )}{\sqrt{U^{n+1}_{k}} + \sqrt{U^{n+1}_{k+1} }} \delta t^{k-1} 
		-{\rm div}\, \bigg( \lm_{w,f} (S^{n+1}) k^*\gr P_{w}^{n+1} \bigg)\\
		&=  \Phi_f \frac{S^{n}}{\delta t^n} + \frac{2C_m}{\delta t^n} 
		\sum_{k=1}^{n}  \frac{\widehat\alpha_m^k (\mathcal{P}(S^k) -\mathcal{P}(S^0) )}{\sqrt{U^{n}_{k} } + \sqrt{U^{n}_{k+1}}}\delta t^{k-1} 
	\end{split}
	\label{disc-w-tipII}
\end{equation}
\begin{equation}
	\begin{split}
		- \Phi_f \frac{S^{n+1}}{\delta t^n} 
		&- \frac{2C_m}{\delta t^n} \sum_{k=1}^{n+1}  \frac{ \widehat\alpha_m^k  (\mathcal{P}(S^k) -\mathcal{P}(S^0) ) }{\sqrt{U^{n+1}_{k}} + \sqrt{U^{n+1}_{k+1} }} \delta t^{k-1}
		- {\rm div}\, \bigg(\lm_{n,f} (S^{n+1}) k^*\gr P_{n}^{n+1} \bigg)\\
		& = - \Phi_f \frac{ S^{n}}{\delta t^n}  - \frac{2C_m}{\delta t^n} 
		\sum_{k=1}^{n}  \frac{ \widehat\alpha_m^k  (\mathcal{P}(S^k) -\mathcal{P}(S^0) ) }{\sqrt{U^{n}_{k} } + \sqrt{U^{n}_{k+1}}}\delta t^{k-1}
	\end{split}
	\label{disc-n-tipII}
\end{equation}
Note that
\begin{align*}
	\widehat\alpha_m^k &=\frac{ \beta_m(\mathcal{P}(S_{\rm max}^{k})) -\beta_m(\mathcal{P}(S_{\rm min}^{k})) }{
		\mathcal{P}(S_{\rm max}^{k})- \mathcal{P}(S_{\rm min}^{k})}.
\end{align*}
where
\begin{align*}
	S_{\rm max}^{k}(x) =\max_{0\leq j\leq k} S^j(x),\quad   S_{\rm min}^{k}(x) =\min_{0\leq j\leq k} S^j(x).
\end{align*}

Using standard finite volume discretization of the two phase system written in the phase formulation (see \cite{EHM-petroleum}) 
we get
\begin{equation}
	\begin{split}
		\Phi_{f,K} \frac{S^{n+1}_K}{\delta t^n} &+ \frac{2C_m}{\delta t^n} 
		\sum_{k=1}^{n+1}  \frac{\widehat\alpha_{m,K}^k (\mathcal{P}(S^k_K) -\mathcal{P}(S^0_K) )}{\sqrt{U^{n+1}_{k, K}} + \sqrt{U^{n+1}_{k+1, K} }} 
		\delta t^{k-1} 
		- \sum_{L\in {\cal N}_K}\tau_{K|L}  k^*_{K|L}  \lm_{w,f,K|L}^{n+1} \delta_{K,L}^{n+1}(P_{w})\\
		&=  \Phi_f \frac{S^{n}_K}{\delta t^n} +  \frac{2C_m}{\delta t^n} 
		\sum_{k=1}^{n}  \frac{\widehat\alpha_{m,K}^k (\mathcal{P}(S^k_K) -\mathcal{P}(S^0_K) )}{\sqrt{U^{n}_{k, K} } + 
			\sqrt{U^{n}_{k+1, K}}}\delta t^{k-1},
	\end{split}
	\label{disc-w-tipII-full}
\end{equation}
\begin{equation}
	\begin{split}
		- \Phi_{f,K} \frac{S^{n+1}_K}{\delta t^n} 
		&- \frac{2C_m}{\delta t^n} \sum_{k=1}^{n+1}  \frac{ \widehat\alpha_{m,K}^k  (\mathcal{P}(S^k_K) -\mathcal{P}(S^0_K) ) }{\sqrt{U^{n+1}_{k, K}}
			+ \sqrt{U^{n+1}_{k+1, K} }} \delta t^{k-1}
		-  \sum_{L\in {\cal N}_K}\tau_{K|L}  k^*_{K|L}  \lm_{n,f,K|L}^{n+1} \delta_{K,L}^{n+1}(P_{n})\\
		& = - \Phi_{f,K} \frac{S^{n}_K}{\delta t^n}  - \frac{2C_m}{\delta t^n} 
		\sum_{k=1}^{n}  \frac{ \widehat\alpha_{m,K}^k  (\mathcal{P}(S^k_K) -\mathcal{P}(S^0_K) ) }{\sqrt{U^{n}_{k, K} }
			+ \sqrt{U^{n}_{k+1}}}\delta t^{k-1, K},
	\end{split}
	\label{disc-n-tipII-full}
\end{equation}
where 
\begin{align}
	U_{k, K}^n = \sum_{l=k}^{n} \widehat\alpha_{m,K}^l \delta t^{l-1}, \label{U-k-n-el}
\end{align}
and 
\begin{align}
	\widehat\alpha_{m,K}^k &=\frac{ \beta_m(\mathcal{P}(S_{\rm max, K}^{k})) -\beta_m(\mathcal{P}(S_{\rm min, K}^{k})) }{
		\mathcal{P}(S_{\rm max, K}^{k})- \mathcal{P}(S_{\rm min, K}^{k})}.
\end{align}
where
\begin{align*}
	S_{\rm max, K}^{k} =\max_{0\leq j\leq k} S^j_K,\quad   S_{\rm min, K}^{k}=\min_{0\leq j\leq k} S^j_K.
\end{align*}


\section{Conclusion}
In this work we have proposed a new, more general way to linearize the imbibition equation (\ref{imb-eqn}) which appears in the definition of the matrix-fracture transfer source terms in the global dual porosity $\delta$-model of incompressible two-phase flow in porous media. After passage to the limit as $\delta\to0$, we analyze the effective matrix-fracture exchange source term obtained by this new linearization and compare it to the effective matrix-fracture exchange source terms obtained previously by a constant linearization in \cite{JPV}. Numerical simulations are provided which show that the matrix-fracture exchange term based on the new linearization procedure gives a better approximation of the exact one than the corresponding exchange term obtained earlier by the authors. Finally, for the effective system in both cases of linearization we provide the discretization schemes by the cell centered finite volume method.

\section{Acknowledgement}
The work of L. Pankratov is supported by Russian Foundation for Basic Research (Grant No. 20-01-00564).
The work of A. Vrba\v{s}ki has been supported by Croatian Science Foundation, project number UIP-2017-05-7249.
Their support is gratefully acknowledged.



\begin{thebibliography}{9999}
	
	\bibitem{ainouz} Ainouz A., Homogenization of a double porosity model
	in deformable media, {\it Electronic Journal of Differential
		Equations}, {\bf 2013}:90 (2013), 1-18.
	
	\bibitem{Chavant-dp}  C. Alboin, J. Jaffr\'{e}, P. Joly, J. Roberts,
	A comparison of methods for calculating the matrix block source term in a double
	porosity model for contaminant transport. Computational Geosciences 6(3) (2002) pp. 523-543.
	
	\bibitem{AMP-B} Amaziane B., Panfilov M., Pankratov L.,
	Homogenized model of two-phase flow with local nonequilibrium in double porosity media, 
	{\it Advances in the Mathematical Physics}, {\bf 2016} (2016), Article ID 3058710, 1-13.
	
	\bibitem{BA-LP-M2AS} B. Amaziane, L. Pankratov, Homogenization of a
	model for water-gas flow through double-porosity media,
	{\it Mathematical Methods in the Applied Sciences}, {\bf 39}:3 (2016), 425-451.
	
	\bibitem{th-layer-M3AS} B. Amaziane, L. Pankratov, A. Piatnitski, Homogenization
	of a single phase flow in a porous medium containing a thin layer,
	{\it Mathematical Models and Methods in Applied Sciences}, {\bf 17}:9 (2007), 1317-1349.
	
	\bibitem{ba-lp-ap} B.~Amaziane, L.~Pankratov and A.~Piatnitski,
	The existence of weak solutions to immiscible compressible two-phase flow
	in porous media: The case of fields with different rock-types,
	{\it Discrete Contin. Dyn. Syst. Ser. B}, {\bf 18}:5 (2013), 1217-1251.
	
	\bibitem{rthin-2} Amaziane B., Pankratov L. , Rybalko V.,
	On the homogenization of some double porosity models with periodic
	thin structures, {\it Applicable Analysis}, {\bf 88} (2009), 1469-1492.
	
	\bibitem{ant-GP-first} S.~N. Antontsev, On the solvability of boundary value problems
	for degenerating equations of two-phase flow, {\it Solid-state Dynamics}, {\bf 10},
	(1972), 28-53 (in Russian).
	
	\bibitem{ant-kaz-mon1990} S.~N. Antontsev, A.~V. Kazhikhov and V.~N. Monakhov,
	{\it Kraevye Zadachi Mekhaniki Neodnorodnykh Zhidkostej},  Nauka, Sibirsk. Otdel.,
	Novosibirsk, 1983 (in Russian); English translation:  {\it Boundary Value Problems in
		Mechanics of Nonhomogeneous Fluids}, North-Holland, Amsterdam, 1990.
	
	\bibitem{Arb-simpl} T. Arbogast,
	{\it A simplified dual-porosity model for two-phase flow},
	in Computational Methods in Water Resources IX, Vol. 2 (Denver, CO, 1992):
	Mathematical Modeling in Water Resources, T.F. Russell, R.E. Ewing, C.A. Brebbia, W.G. Gray, and
	G.F. Pindar, eds., Comput. Mech., Southampton, U.K., 1992, pp. 419-426.
	
	\bibitem{Arbogas-90} T. Arbogast, J. Douglas, Jr.,  U. Hornung, Derivation of the double porosity
	model of single phase flow via homogenization theory, Siam J. Math. Anal., 21 (1990), pp. 823-836.
	
	\bibitem{ADH91} T. Arbogast, J. Douglas, U. Hornung,
	{\it Modeling of naturally fractured reservoirs by formal homogenization techniques},
	in Frontiers in Pure and Applied Mathematics, R. Dautray, ed.,
	North-Holland, Amsterdam, 1991, pp. 1--19.
	
	\bibitem{bp} N. S. Bakhvalov, G.P. Panasenko, {\it Averaging processes in
		periodic media}, Nauka, Moscow, 1984; English transl., Kluwer, Dordrecht, 1989.
	
	\bibitem{BZK-1960}  G. Barenblatt, I. Zheltov, and I. Kochina, Basic concepts in the theory of seepage
	of homogeneous liquids in the fractured rock, J. Appl. Math. Mech., 24 (1960), pp. 1286-1303.
	
	\bibitem{Bear-1993}    Bear, J., C. F. Tsang, and G. de Marsily (1993), Flow and Contaminant
	Transport in Fractured Rock, Academic, San Diego, Calif.
	
	\bibitem{bcp} Bourgeat, A., Chechkin, G. A. and Piatnitski, A., 2003,
	Singular double porosity model. {\itshape Appl. Anal.}, {\bfseries 82}, 103-116.
	
	\bibitem{bgpp} Bourgeat, A., Goncharenko, M., Panfilov, M. and Pankratov, L.,
	1999, A general double porosity model.
	{\itshape C. R. Acad. Sci. Paris, S\'erie IIb},  {\bfseries 327}, 1245-1250.
	
	\bibitem{bmp} Bourgeat, A., Mikelic, A. and Piatnitski, A., 1998,
	Mod\`ele de double porosit\'e al\'eatoire.
	{\itshape C. R. Acad. Sci. Paris, S\'erie 1}, {\bfseries 327}, 99-104.
	
	\bibitem{BLM} A. Bourgeat, S. Luckhaus, A. Mikeli\'c,
	{\it Convergence of the homogenization process for a double-porosity
		model of immiscible two-phase flow}, SIAM J. Math. Anal., {\bf 27} (6) (1996), pp.~1520-1543.
	
	\bibitem{brad} Braides, A., Chiad\`o Piat, V. and Piatnitski, A., 2004,
	A variational approach to double-porosity problems.
	{\itshape Asymptotic Anal.}, {\bfseries 39}, 281-308.
	
	\bibitem{GC-JJ} G.~Chavent and J.~Jaffr\'e, {\it Mathematical Models and Finite
		Elements for Reservoir Simulation}, North-Holland, Amsterdam, 1986.
	
	\bibitem{ZC-GH-YM-06} Z.~Chen, G.~Huan and Y.~Ma,
	{\it Computational Methods for Multiphase Flows in Porous Media}, SIAM, Philadelphia, 2006.
	
	\bibitem{cat} Choquet, C., 2004,
	Derivation of the double porosity model of a compressible miscible displacement
	in naturally fractured reservoirs. {\itshape Appl. Anal.}, {\bfseries 83}, 477-499.
	
	\bibitem{cior} Cioranescu D, Saint Jean Paulin J.
	{\it Homogenization of Reticulated Structures}. Applied Mathematical Sciences, Vol. 136. Springer-Verlag: New York-Berlin-Heidelberg. 1999.
	
	\bibitem{EGH-book} R. Eymard, T. Gallou\"et, R. Herbin, Finite volume methods, in
	Handbook of numerical analysis, Vol. VII, pages 713-1020, North-Holland, Amsterdam, 2000.
	
	\bibitem{EHM-petroleum} R. Eymard, R. Herbin, A. Michel,
	Mathematical study of a petroleum-engineering scheme,
	ESAIM: M2AN, vol 37,  no 6, 2003, pp 937-972.
	
	\bibitem{CG-MS2} C.~Galusinski and M.~Saad, Water-gas flow in porous media,
	{\it Discrete Contin. Dyn. Syst. Ser. B}, {\bf 9} (2008), 281-308.
	
	\bibitem{CG-MS4} C.~Galusinski and M.~Saad,
	Weak solutions for immiscible compressible multifluid flows in porous media,
	{\it C. R. Acad. Sci. Paris, S\'er. I}, {\bf 347} (2009), 249-254.
	
	\bibitem{Hornung-97}  U. Hornung, ed., Homogenization and Porous Media,
	no. 6 in Interdisciplinary Applied Mathematics, Springer-Verlag, New York, 1997.
	
	\bibitem{Firoo-2005} Hoteit, H.  Firoozabadi, A.
	An efficient numerical model for incompressible two-phase flow in fractured media,
	Advances in Water Resources 31 (2008) 891-905.
	
	\bibitem{JPV} M.~Jurak, L.~Pankratov, A.~Vrba\v{s}ki, A fully homogenized model for
	incompressible two-phase flow in double porosity media, {\it Applicable Analysis},
	{\bf 95}:10 (2016), 2280-2299.
	
	\bibitem{khal-saad} Z.~Khalil and M.~Saad, Solutions to a model for compressible
	immiscible two phase flow in porous media, {\it Electronic Journal of Differential Equations},
	{\bf 122} (2010), 1-33.
	
	\bibitem{khal-saad2} Z.~Khalil and M.~Saad, On a fully nonlinear degenerate
	parabolic system modeling immiscible gas-water displacement in porous media,
	{\it Nonlinear Analysis: Real World Applications}, {\bf 12} (2011), 1591-1615.
	
	\bibitem{ak-lp-AA} A.~Konyukhov, L.~Pankratov, Upscaling of an immiscible non-equilibrium
	two-phase flow in double porosity media, {\it Applicable Analysis}, {\bf 95} (2016),
	2300-2322.
	
	\bibitem{Linss} Torsten Linss, {\it Layer-Adapted Meshes for
		Reaction-Convection-Diffusion Problems}, Springer, 2010.
	
	\bibitem{kh} Marchenko, V. A. and Khruslov, E. Ya., 2006,
	{\itshape Homogenization of Partial Differential Equations},
	(Boston: Birkh\"auser).
	
	\bibitem{Nelson-2001} Nelson, R. A. (2001). Geologic Analysis of Naturally Fractured
	Reservoirs. Gulf Professional Publishing, Boston, Massachusetts.
	
	\bibitem{panf}  Panfilov, M., 2000,
	{\itshape Macroscale Models of Flow Through Highly Heterogeneous Porous Media},
	(Dordrecht-Boston-London: Kluwer Academic Publishers).
	
	\bibitem{pr} Pankratov, L. and Rybalko, V., 2003,
	Asymptotic analysis of a double porosity model with thin fissures.
	{\itshape Mat. Sbornik},  {\bfseries 194}, 121-146.
	
	\bibitem{RAG} R. Raghavan, E. Ozkan, {\it A Method for Computing Unsteady Flows in Porous Media},
	Pitman Research Notes in Mathematics, 318, Longman Scientific and Technical, 1994.
	
	\bibitem{Bastian-06} V. Reichenberger, H. Jakobs, P. Bastien, R. Helmig,
	A mixed-dimensional finite volume method for multiphase flow in fractured porous
	media, Adv. Water Resources 29 (7) (2006) 1020-1036.
	
	\bibitem{san} Sandrakov, G. V., 1999,
	Homogenization of parabolic equations with contrasting coefficients.
	{\itshape Izv. Math.}, {\bfseries 63}, 1015-1061.
	
	\bibitem{deSwaan} de Swaan, A. (1978), Theory of waterflooding in fractured reservoirs,
	SPE J., 18, 117-122.
	
	\bibitem{} C. J. van Duijn, J. Molenaar, M. de Neef, The effect of capillary forces on
	immiscible two-phase flow in heteregeneous porous media, Transport in Porous Media 21 (1995) 71-93.
	
	\bibitem{WarrenRoot} J. Warren and P. Root, The behavior of naturally fractured reservoirs, Soc. Pet.
	Eng., J. 3 (1963), pp. 245-255.
	
	\bibitem {Yeh06} Li--Ming Yeh,
	Homogenization of two-phase flow in fractured media,
	Math. Models Methods Appl. Sci., {\bf 16} (10) (2006), pp.~1627-1651.
	
	
	
\end{thebibliography}
\end{document}